\documentclass{amsart}

\usepackage{amsthm}
\usepackage[leqno]{amsmath}
\usepackage[all]{xy} \SelectTips{eu}{}
\usepackage{latexsym,amsfonts,amssymb}
\usepackage{hyperref}

\newcommand{\numberseries}{\bfseries}   
\newlength{\thmtopspace}                
\newlength{\thmbotspace}                
\newlength{\thmheadspace}               
\newlength{\thmindent}                  
\newtheoremstyle{bfupright head,slanted body}
                {\thmtopspace}{\thmbotspace}
                {\slshape}{\thmindent}{\bfseries}{.}{\thmheadspace}
                {{\numberseries \thmnumber{#2\;}}\thmnote{#3}}

\newtheoremstyle{bfupright head,upright body}
                {\thmtopspace}{\thmbotspace}
                {\upshape}{\thmindent}{\bfseries}{.}{\thmheadspace}
                {{\numberseries \thmnumber{#2\;}}\thmnote{#3}}

\newtheoremstyle{fixed bf head,slanted body}
                {\thmtopspace}{\thmbotspace}{\slshape}
                {\thmindent}{\bfseries}{.}{\thmheadspace}
                {{\numberseries \thmnumber{#2\;}}\thmname{#1}\thmnote{ (#3)}}

\newtheoremstyle{fixed bf head,upright body}
                {\thmtopspace}{\thmbotspace}{\upshape}
                {\thmindent}{\bfseries}{.}{\thmheadspace}
                {{\numberseries \thmnumber{#2\;}}\thmname{#1}\thmnote{ (#3)}}

\newtheoremstyle{numbered paragraph}
                {\thmtopspace}{\thmbotspace}{\upshape}
                {\thmindent}{\upshape}{}{\thmheadspace}
                {{\numberseries \thmnumber{#2.}}}

\theoremstyle{bfupright head,slanted body}
\newtheorem{res}{}[section]             \newtheorem*{res*}{}

\theoremstyle{bfupright head,upright body}
\newtheorem{bfhpg}[res]{}               \newtheorem*{bfhpg*}{}

\theoremstyle{fixed bf head,slanted body}
\newtheorem{thm}[res]{Theorem}          \newtheorem*{thm*}{Theorem}
\newtheorem{prp}[res]{Proposition}      \newtheorem*{prp*}{Proposition}
\newtheorem{cor}[res]{Corollary}        \newtheorem*{cor*}{Corollary}
\newtheorem{lem}[res]{Lemma}            \newtheorem*{lem*}{Lemma}

\theoremstyle{fixed bf head,upright body}
\newtheorem{dfn}[res]{Definition}       \newtheorem*{dfn*}{Definition}
\newtheorem{rmk}[res]{Remark}           \newtheorem*{rmk*}{Remark}

\theoremstyle{numbered paragraph}
\newtheorem{ipg}[res]{}

\newlength{\thmlistleft}        
\newlength{\thmlistright}       
\newlength{\thmlistpartopsep}   
\newlength{\thmlisttopsep}      
\newlength{\thmlistparsep}      
\newlength{\thmlistitemsep}     

\newcounter{eqc}
\newenvironment{eqc}{\begin{list}{\upshape (\textit{\roman{eqc}})}%
    {\usecounter{eqc}%
      \setlength{\leftmargin}{\thmlistleft}%
      \setlength{\labelwidth}{\thmlistleft}%
      \setlength{\rightmargin}{\thmlistright}%
      \setlength{\partopsep}{\thmlistpartopsep}%
      \setlength{\topsep}{\thmlisttopsep}%
      \setlength{\parsep}{\thmlistparsep}%
      \setlength{\itemsep}{\thmlistitemsep}}}%
  {\end{list}}%

\newcounter{prt}
\newenvironment{prt}{\begin{list}{\upshape (\alph{prt})}%
    {\usecounter{prt}%
      \setlength{\leftmargin}{\thmlistleft}%
      \setlength{\labelwidth}{\thmlistleft}%
      \setlength{\rightmargin}{\thmlistright}%
      \setlength{\partopsep}{\thmlistpartopsep}%
      \setlength{\topsep}{\thmlisttopsep}%
      \setlength{\parsep}{\thmlistparsep}%
      \setlength{\itemsep}{\thmlistitemsep}}}%
  {\end{list}}%

\newcounter{rqm}
\newenvironment{rqm}{\begin{list}{\upshape (\arabic{rqm})}%
    {\usecounter{rqm}%
      \setlength{\leftmargin}{\thmlistleft}%
      \setlength{\labelwidth}{\thmlistleft}%
      \setlength{\rightmargin}{\thmlistright}%
      \setlength{\partopsep}{\thmlistpartopsep}%
      \setlength{\topsep}{\thmlisttopsep}%
      \setlength{\parsep}{\thmlistparsep}%
      \setlength{\itemsep}{\thmlistitemsep}}}%
  {\end{list}}%

\newenvironment{prf*}[1][Proof]{%
  \begin{proof}[\bf #1]
    \setcounter{equation}{0}
    }
  {\end{proof}
}

\newcommand{\pgref}[1]{\ref{#1}}
\newcommand{\thmref}[2][Theorem~]{#1\pgref{thm:#2}}
\newcommand{\corref}[2][Corollary~]{#1\pgref{cor:#2}}
\newcommand{\prpref}[2][Proposition~]{#1\pgref{prp:#2}}
\newcommand{\lemref}[2][Lemma~]{#1\pgref{lem:#2}}
\newcommand{\dfnref}[2][Definition~]{#1\pgref{dfn:#2}}
\newcommand{\rmkref}[2][Remark~]{#1\pgref{rmk:#2}}
\renewcommand{\eqref}[1]{(\pgref{eq:#1})}

\newcommand{\thmcite}[2][?]{\cite[thm.~#1]{#2}}
\newcommand{\prpcite}[2][?]{\cite[prop.~#1]{#2}}
\newcommand{\corcite}[2][?]{\cite[cor.~#1]{#2}}
\newcommand{\lemcite}[2][?]{\cite[lem.~#1]{#2}}
\newcommand{\seccite}[2][?]{\cite[sec.~#1]{#2}}
\newcommand{\dfncite}[2][?]{\cite[def.~#1]{#2}}
\newcommand{\rmkcite}[2][?]{\cite[rmk.~#1]{#2}}

\newcommand{\eqclbl}[1]{{\upshape(\textit{#1})}}
\newcommand{\proofofimp}[3][:]{\mbox{\eqclbl{#2}$\!\implies\!$\eqclbl{#3}#1}}

\numberwithin{equation}{res}

\def\urltilda{\kern -.15em\lower .7ex\hbox{\~{}}\kern .04em}
\newcommand{\set}[2][\mspace{1mu}]{\{#1 #2 #1\}}
\newcommand{\setof}[3][\mspace{1mu}]{\{#1#2 \mid #3#1\}}
\newcommand{\ZZ}{\mathbb{Z}}
\newcommand{\deq}{\:=\:}

\newcommand{\dis}{\:\is\:}

\renewcommand{\a}{\alpha}
\newcommand{\f}{\varphi}
\newcommand{\is}{\cong}
\newcommand{\qis}{\simeq}
\renewcommand{\le}{\leqslant}
\renewcommand{\ge}{\geqslant}

\newcommand{\lra}{\longrightarrow}
\newcommand{\xra}[2][]{\xrightarrow[#1]{\;#2\;}}
\newcommand{\qra}{\xra{\qis}}
\newcommand{\QQ}{\mathbb{Q}}
\newcommand{\Rop}{R^\circ}
\newcommand{\Sop}{S^\circ}
\newcommand{\mapdef}[4][\rightarrow]{\nobreak{#2\colon #3 #1 #4}}

\renewcommand{\Im}[1]{\nobreak{\operatorname{Im}#1}}
\newcommand{\Ker}[1]{\nobreak{\operatorname{Ker}#1}}
\newcommand{\dif}[2][]{{\partial}^{#2}_{#1}}
\newcommand{\Bo}[2][]{\operatorname{B}_{#1}(#2)}
\newcommand{\Cy}[2][]{\operatorname{Z}_{#1}(#2)}
\newcommand{\Co}[2][]{\operatorname{C}_{#1}(#2)}
\renewcommand{\H}[2][]{\operatorname{H}_{#1}(#2)}
\newcommand{\HH}[2][]{\operatorname{H}^{#1}(#2)}
\newcommand{\Shift}[2][]{\mathsf{\Sigma}^{#1}{#2}}
\newcommand{\Thb}[2]{#2_{{\scriptscriptstyle\ge}#1}}
\newcommand{\fd}[2][R]{\operatorname{fd}_{#1}#2}
\newcommand{\pd}[2][R]{\operatorname{proj.\mspace{-2mu}dim}_{#1}#2}
\newcommand{\Gfd}[2][R]{\operatorname{Gfd}_{#1}#2}
\newcommand{\Gid}[2][R]{\operatorname{Gid}_{#1}#2}
\newcommand{\Hom}[3][R]{\operatorname{Hom}_{#1}(#2,#3)}
\newcommand{\RHom}[3][R]{\operatorname{\mathbf{R}Hom}_{#1}(#2,#3)}
\newcommand{\Ext}[4][R]{\operatorname{Ext}_{#1}^{#2}(#3,#4)}
\newcommand{\tp}[3][R]{\nobreak{#2\otimes_{#1}#3}}
\newcommand{\tpp}[3][R]{(\tp[#1]{#2}{#3})}
\newcommand{\Ltp}[3][R]{\nobreak{#2\otimes_{#1}^{\mathbf{L}}#3}}
\newcommand{\Ltpp}[3][R]{(\Ltp[#1]{#2}{#3})}
\newcommand{\Tor}[4][R]{\operatorname{Tor}^{#1}_{#2}(#3,#4)}
\newcommand{\xycomma}[1][,]{\rlap{\;#1}}
\newcommand{\ddif}[2][]{{\partial}_{#2}^{#1}}
\newcommand{\dimR}[1][R]{\operatorname{dim}#1}
\newcommand{\splfR}[1][R]{\operatorname{splf}#1}
\newcommand{\HBo}[2][]{\operatorname{B}^{#1}(#2)}
\newcommand{\HCy}[2][]{\operatorname{Z}^{#1}(#2)}
\newcommand{\HCo}[2][]{\operatorname{C}^{#1}(#2)}

   \def\soft#1{\leavevmode\setbox0=\hbox{h}\dimen7=\ht0\advance
    \dimen7 by-1ex\relax\if t#1\relax\rlap{\raise.6\dimen7
    \hbox{\kern.3ex\char'47}}#1\relax\else\if T#1\relax
    \rlap{\raise.5\dimen7\hbox{\kern1.3ex\char'47}}#1\relax
    \else\if d#1\relax\rlap{\raise.5\dimen7\hbox{\kern.9ex
    \char'47}}#1\relax\else\if D#1\relax\rlap{\raise.5\dimen7
    \hbox{\kern1.4ex\char'47}}#1\relax\else\if l#1\relax
    \rlap{\raise.5\dimen7\hbox{\kern.4ex\char'47}}#1\relax
    \else\if L#1\relax\rlap{\raise.5\dimen7\hbox{\kern.7ex
    \char'47}}#1\relax\else\message{accent \string\soft
    \space #1 not defined!}#1\relax\fi\fi\fi\fi\fi\fi}

\begin{document}

\title[Gorenstein dimensions of unbounded complexes and change of
base]%
{Gorenstein dimensions of unbounded complexes and change of base (With
  an appendix\\ by Driss Bennis)}

\author{Lars Winther Christensen}

\address{L.W.C. Texas Tech University, Lubbock, TX 79409, U.S.A.}

\email{lars.w.christensen@ttu.edu}

\urladdr{http://www.math.ttu.edu/\urltilda lchriste}

\author{Fatih K\"{o}ksal}

\address{F.K. Lewis University, Romeoville, IL 60446, U.S.A.}

\email{koksalfa@lewisu.edu}

\author{Li Liang}

\address{L.L. Lanzhou Jiaotong University, Lanzhou 730070, China}

\email{lliangnju@gmail.com}

\thanks{This research was partly supported by NSA grant H98230-14-0140
  (L.W.C.), by NSFC grant 11301240 and SRF for ROCS, SEM (L.L.). Part
  of the work was done in the summer of 2015 when L.W.C.\ and L.L.\
  visited Nanjing University supported by NSFC grant 11371187 (PI
  N.Q.~Ding).}

\date{12 May 2016}



\begin{abstract}
  For a commutative ring $R$ and a faithfully flat $R$-algebra $S$ we
  prove, under mild extra assumptions, that an $R$-module $M$ is
  Gorenstein flat if and only if the left $S$-module $\tp{S}{M}$ is
  Gorenstein flat, and that an $R$-module $N$ is Gorenstein injective
  if and only if it is cotorsion and the left $S$-module $\Hom{S}{N}$
  is Gorenstein injective.  We apply these results to the study
  of Gorenstein homological dimensions of unbounded complexes. In
  particular, we prove two theorems on stability of these dimensions
  under faithfully flat (co-)base change.
\end{abstract}

\maketitle

\thispagestyle{empty}

\section{Introduction}

\noindent Auslander and Bridger's \cite{MAsMBr69} notion of
G-dimension for finitely generated modules over noetherian rings was
generalized and dualized by Enochs and collaborators, who introduced
the Gorenstein projective, Gorenstein injective, and Gorenstein flat
dimension of modules over associative rings. In the treatment by
Christensen, Frankild, and Holm \cite{CFH-06}, these invariants were
considered for complexes with bounded homology. It is, however,
possible to define Gorenstein projective dimension for unbounded
complexes: This was done by Veliche \cite{OVl06}, and the dual case of
Gorenstein injective dimension was treated by Asadollahi and Salarian
\cite{JAsSSl06a}.

Inspired by \cite{OVl06}, we propose a definition of Gorenstein flat
dimension for unbounded complexes; it coincides with the one
introduced by Iacob~\cite{AIc09} whenever the latter is defined. Our
main results are two theorems on stability of Gorenstein homological
dimensions of complexes under faithfully flat change of base.

\begin{thm}
  \label{thm:A}
  Let $R$ be a commutative coherent ring and let $S$ be a faithfully
  flat $R$-algebra that is left GF-closed. For every $R$-complex $M$
  there is an equality
  \begin{equation*}
    \Gfd{M} \deq \Gfd[S]{\Ltpp{S}{M}}\:.
  \end{equation*}
  In particular, an $R$-module $M$ is Gorenstein flat if and only if
  the $S$-module $\tp{S}{M}$ is Gorenstein flat.
\end{thm}
\noindent
This result is part of \thmref{transfer}; the condition that $S$ is
left GF-closed is discussed in \pgref{GF-closed}; it is satisfied if
$S$ is right coherent.

\begin{thm}
  \label{thm:B}
  Let $R$ be a commutative noetherian ring with $\splfR <\infty$,
  i.e.\ every flat $R$-module has finite projective dimension, and let
  $S$ be a faithfully flat $R$-algebra. For every $R$-complex $N$ with
  $\HH[i]{N}=0$ for all $i \gg 0$ there is an equality
  \begin{equation*}
    \Gid{N} \deq \Gid[S]{\RHom{S}{N}}\:.
  \end{equation*}
  In particular, an $R$-module $N$ is Gorenstein injective if and only
  if the $S$-module $\Hom{S}{N}$ is Gorenstein injective and
  $\Ext{i}{S}{N}=0$ holds for all $i >0$.
\end{thm}
\noindent
This result is part of \thmref{3}. The condition $\splfR < \infty$ is
discussed in \pgref{splf}; it is satisfied if $R$ has finite Krull
dimension or cardinality at most $\aleph_n$ for some integer $n\ge
0$. We notice that the condition implies vanishing of $\Ext{i}{S}{N}$
for~$i \gg 0$.

A recent result of \v{S}{\soft{t}}ov{\'{\i}}{\v{c}}ek \cite{JSt}
implies that every Gorenstein injective module is cotorsion. It allows
us to prove the next result, which appears as \thmref{2}.
\begin{thm}
  \label{thm:C}
  Let $R$ be commutative noetherian and let $S$ be a faithfully flat
  $R$-algebra. An $R$-module $N$ is Gorenstein injective if and only
  if it is cotorsion and the $S$-module $\Hom{S}{N}$ is Gorenstein
  injective.
\end{thm}
\noindent This result compares to the statement about modules in
\thmref{B}.

By significantly relaxing the conditions on the rings, the results of
this paper improve results obtained by Christensen and
Holm~\cite{LWCHHl09}, by Christensen and
Sather-Wagstaff~\cite{LWCSSW10}, and by Liu and
Ren~\cite{LZhRWe14}. Details pertaining to \thmref[Theorems~]{A} and
\thmref[]{B} are given in \rmkref[Remarks~]{1} and \rmkref[]{2}; the
trend is that the rings in \cite{LWCSSW10,LZhRWe14} are assumed to be
commutative noetherian and, more often than not, of finite Krull
dimension.

The paper is organized as follows: In section 2 we set the notation
and recall some background material. Sections 3--4 focus on the
Gorenstein injective dimension, and Sections 5--6 deal with the
Gorenstein flat dimension. Section 7 has some closing remarks and,
finally, an appendix by Bennis answers a question raised in an earlier
version of this paper.

\section{Complexes}

\noindent
Let $R$ be a ring with identity. We consider only unitary $R$-modules,
and we employ the convention that $R$ acts on the left. That is, an
$R$-module is a left $R$-module, and right $R$-modules are treated as
modules over the opposite ring, denoted~$\Rop$.

\begin{bfhpg}[\bf Complexes]
  Complexes of $R$-modules, $R$-complexes for short, is our object of
  study. Let $M$ be an $R$-complex. With homological grading, $M$ has
  the form
  \begin{equation*}
    \cdots \lra M_{i+1} \xra{\dif[i+1]{M}} M_i \xra{\dif[i]{M}}
    M_{i-1} \lra \cdots\: ;
  \end{equation*}
  one switches to cohomological grading by setting $M^i = M_{-i}$ and
  $\ddif[i]{M} = \dif[-i]{M}$ to get
  \begin{equation*}
    \cdots \lra M^{i-1} \xra{\ddif[i-1]{M}} M^i \xra{\ddif[i]{M}}
    M^{i+1} \lra \cdots\:.
  \end{equation*}

  For $n\in\ZZ$ the symbol $M_{\ge n}$ denotes the quotient complex of
  $M$ with $(M_{\ge n})_i=M_i$ for $i\ge n$ and $(M_{\ge n})_i=0$ for
  $i < n$.

  The subcomplexes $\Bo{M}$ and $\Cy{M}$ of boundaries and cycles, the
  quotient complex $\Co{M}$ of cokernels, and the subquotient complex
  $\H{M}$ of homology all have zero differentials; their modules are
  given by
  \begin{alignat*}{3}
    \Bo[i]{M} &= \Im{\dif[i+1]{M}}& &=& \Im{\ddif[-i-1]{M}} &= \HBo[-i]{M}\\
    \Cy[i]{M} &= \Ker{\dif[i]{M}}& &=& \Ker{\ddif[-i]{M}} &= \HCy[-i]{M}\\
    \Co[i]{M} &= M_i/\Bo[i]{M}& &=& M^{-i}/\HBo[-i]{M} &= \HCo[-i]{M}\\
    \H[i]{M} &= \Cy[i]{M}/\Bo[i]{M}& &=&\;\HCy[-i]{M}/\HBo[-i]{M} &=
    \HH[-i]{M}
  \end{alignat*}

  The complex $M$ is called \emph{acyclic} if $\H{M}=0$, i.e.\
  $\H[i]{M}=0$ holds for all $i\in\ZZ$.
\end{bfhpg}

\begin{bfhpg}[\bf Morphisms]
  A morphism $\mapdef{\a}{M}{N}$ of $R$-complexes is a family
  of $R$-linear maps $\set{\mapdef{\a_i}{M_i}{N_i}}_{i\in\ZZ}$  that
  commute with the differentials. That is, for all $i\in\ZZ$ one has
  $\a_{i-1}\dif[i]{M} = \dif[i]{N}\a_i$. A morphism
  $\mapdef{\a}{M}{N}$ maps boundaries to boundaries and cycles to
  cycles, so it induces a morphism $\mapdef{\H{\a}}{\H{M}}{\H{N}}$. If
  the induced morphism $\H{\a}$ is bijective, then $\a$ is called a
  \emph{quasi-isomorphism.} Such morphisms are marked by the symbol
  `$\qis$'.
\end{bfhpg}

\begin{bfhpg}[\bf Resolutions]
  \label{res}
  An $R$-complex $I$ is called \emph{semi-injective} if each module
  $I^i$ is injective, and the functor $\Hom{-}{I}$ preserves
  quasi-isomorphisms.  Dually, an $R$-complex $P$ is called
  \emph{semi-projective} if each module $P_i$ is projective, and the
  functor $\Hom{P}{-}$ preserves quasi-isomorphisms. An $R$-complex
  $F$ is called \emph{semi-flat} if each module $F_i$ is flat, and the
  functor $\tp{-}{F}$ preserves quasi-isomorphisms. Every
  semi-projective complex is semi-flat. For the following facts see
  Avramov and Foxby~\cite{LLAHBF91}.
  \begin{prt}
  \item[$\bullet$] Every $R$-complex $M$ has a semi-projective
    resolution in the following sense: There is a quasi-isomorphism
    $\mapdef{\pi}{P}{M}$, where $P$ is a semi-projective complex;
    moreover, $\pi$ can be chosen surjective.

  \item[$\bullet$] Every $R$-complex $M$ has a semi-injective
    resolution in the following sense: There is a quasi-isomorphism
    $\mapdef{\iota}{M}{I}$, where $I$ is semi-injective; moreover,
    $\iota$ can be chosen injective.
  \end{prt}
\end{bfhpg}

\begin{bfhpg}[\bf The derived category]
  In the derived category over $R$, the objects are $R$-complexes, and
  the morphisms are equivalence classes of diagrams \mbox{$\bullet
    \xleftarrow{\;\qis\;} \bullet \lra \bullet$} of morphisms of
  $R$-complexes. The isomorphisms in the derived category are classes
  represented by diagrams with two quasi-isomorphisms; and they are
  also marked by the symbol `$\qis$'.

  The derived tensor product, $\Ltp{-}{-}$, and the derived Hom,
  $\RHom{-}{-}$, are functors on the derived category; their values on
  given $R$-complexes are computed by means of the resolutions
  described above.  As for modules one sets
  \begin{equation*}
    \Ext{i}{M}{N} \deq \HH[i]{\RHom{M}{N}}
  \end{equation*}
  for $R$-complexes $M$ and $N$ and $i\in\ZZ$.
\end{bfhpg}

\section{Gorenstein injective modules and cotorsion}

\noindent
For definitions and standard results on Gorenstein homological
dimensions our references are Holm's~\cite{HHl04a} and the monograph
\cite{rha} by Enochs and Jenda. We start by recalling the definition
\dfncite[10.1.1]{rha} of a Gorenstein injective module.

\begin{ipg}
  A complex $U$ of injective $R$-modules is called \emph{totally
    acyclic} if it is acyclic and the complex $\Hom{I}{U}$ is acyclic
  for every injective $R$-module $I$.

  An $R$-module $G$ is called \emph{Gorenstein injective} (for short,
  G-injective) if there exists a totally acyclic complex $U$ of
  injective $R$-modules with $\HCy[0]{U} \is G$.
\end{ipg}

\begin{ipg}
  \label{g-injective}
  Every $R$-module has an injective resolution, so to prove that a
  module $N$ is Gorenstein injective it suffices to verify the
  following:
  \begin{rqm}
  \item $\Ext{i}{I}{N}=0$ holds for all $i>0$ and every injective
    $R$-module $I$.
  \item $N$ has a proper left injective resolution. That is, there
    exists an acyclic complex of $R$-modules $U^+ = \cdots \to U^{-2}
    \to U^{-1} \to N \to 0$ with each $U^i$ injective, such that
    $\Hom{I}{U^+}$ is acyclic for every injective $R$-module $I$.
  \end{rqm}
\end{ipg}

\begin{ipg}
  \label{Jst}
  An $R$-module $X$ is called \emph{cotorsion} if one has
  $\Ext{1}{F}{X}=0$ (equivalently, $\Ext{>0}{F}{X}=0$) for every flat
  $R$-module $F$.

  \v{S}{\soft{t}}ov{\'{\i}}{\v{c}}ek \corcite[5,9]{JSt} shows that for
  every module $N$ that has a left injective resolution, one has
  $\Ext{i}{F}{N}=0$ for all $i>0$ and every $R$-module $F$ of finite
  flat dimension. In particular, every Gorenstein injective module is
  cotorsion.
\end{ipg}

This provides for the following improvement of \cite[Ascent table
II.(h)]{LWCHHl09}

\begin{lem}
  \label{lem:giascent}
  Let $R\to S$ be a ring homomorphism such that $\fd{S}$ and
  $\fd[\Rop\!]{S}$ are finite. For a Gorenstein injective $R$-module
  $G$, the $S$-module $\Hom{S}{G}$ is Gorenstein injective.
\end{lem}

\begin{prf*}
  Let $U$ be a totally acyclic complex of injective $R$-modules with
  $\HCy[0]{U} \is G$. The complex $\Hom{S}{U}$ consists of injective
  $S$-modules, and it is acyclic as $\fd{S}$ is finite; see
  \pgref{Jst}. One has $\HCy[0]{\Hom{S}{U}} \is \Hom{S}{G}$, so it
  remains to show that the complex
  \begin{equation*}
    \Hom[S]{J}{\Hom{S}{U}} \dis \Hom{J}{U}
  \end{equation*}
  is acyclic for every injective $S$-module $J$. As $\fd[\Rop\!]{S}$
  is finite, an injective $S$-module has finite injective dimension
  over $R$; see \corcite[4.2]{LLAHBF91}. Acyclicity of $\Hom{J}{U}$
  now follows from \lemcite[2.2]{CFH-06}.
\end{prf*}

\begin{ipg}
  Let $N$ be an $R$-complex. A \emph{complete injective resolution} of
  $N$ is a diagram
  \begin{equation*}
    N \xra{\iota} I \xra{\upsilon} U
  \end{equation*}
  where $\iota$ is a semi-injective resolution, $\upsilon^i$ is
  bijective for all $i \gg 0$, and $U$ is a totally acyclic complex of
  injective $R$-modules.

  The \emph{Gorenstein injective dimension} of an $R$-complex $N$ is
  defined \dfncite[2.2]{JAsSSl06a} as
  \begin{equation*}
    \Gid{N} \deq \inf\left\{ g \in \ZZ \:
      \left|
        \begin{array}{c}
          N \xra{\iota} I \xra{\upsilon} U\\
          \text{is a complete injective resolution with}\\
          \upsilon^i\colon I^i \to U^i \text{ bijective for all } i\ge g
        \end{array}
      \right.
    \right\}.
  \end{equation*}
\end{ipg}

The next result is dual to (parts of) \thmcite[3.4]{OVl06}, and the
proof is omitted.\footnote{ The statement \thmcite[2.3]{JAsSSl06a} is
  also meant to be the dual of \thmcite[3.4]{OVl06}, but it has an
  unfortunate typographical error, reproduced in
  \lemcite[1]{LZhRWe14}: In parts \eqclbl{ii} and \eqclbl{iii} the
  kernel has been replaced by the cokernel in the same degree.}
\pagebreak
\begin{prp}
  \label{prp:ovl}
  Let $N$ be an $R$-complex and $g$ be an integer. The following
  conditions are equivalent.
  \begin{eqc}
  \item $\Gid{N} \le g$.
  \item $\HH[i]{N} = 0$ holds for all $i > g$, and there exists a
    semi-injective resolution $N \qra I$ such that the module
    $\HCy[g]{I}$ is Gorenstein injective.
  \item $\HH[i]{N} = 0$ holds for all $i > g$, and for every
    semi-injective resolution $N \qra I$ the module $\HCy[g]{I}$ is
    Gorenstein injective.\qed
  \end{eqc}
\end{prp}

The next result improves \corcite[9]{LZhRWe14} by removing assumptions
that $R$ and $S$ should be commutative noetherian with $\dimR<\infty$.

\begin{prp}
  \label{prp:gidascent}
  Let $R\to S$ be a ring homomorphism such that $\fd{S}$ and
  $\fd[\Rop\!]{S}$ are finite. For every $R$-complex $N$ one has
  $\Gid[S]{\RHom{S}{N}} \le \Gid{N}$.
\end{prp}

\begin{prf*}
  We may assume that $\H{N}$ is non-zero, otherwise there is nothing
  to prove. Let $N \qra I$ be a semi-injective resolution; the
  $S$-complex $\Hom{S}{I}$ is semi-injective and isomorphic to
  $\RHom{S}{N}$ in the derived category. If $\Gid{N}$ is finite, say
  $g$, then by \prpref{ovl} one has $\HH[i]{I}=0$ for $i > g$ and $Z =
  \HCy[g]{I}$ is a G-injective $R$-module. Now \pgref{Jst} yields
  $\HH[g+n]{\Hom{S}{I}} = \Ext{n}{S}{Z} = 0$ for all $n>0$. Moreover,
  the $S$-module $\HCy[g]{\Hom{S}{I}} \is \Hom{S}{\HCy[g]{I}}$ is
  G-injective by \lemref{giascent}, so $\Gid[S]{\RHom{S}{N}} \le g$
  holds by \prpref{ovl}.
\end{prf*}

\section{Faithfully flat co-base change}

\noindent Throughout this section, $R$ is a commutative ring and $S$
is an $R$-algebra. We are primarily concerned with the following
setup.

\begin{ipg}
  \label{ff}
  Let $S$ be a faithfully flat $R$-algebra; there is then a pure exact
  sequence
  \begin{equation}
    \label{eq:seq}
    0 \lra R \lra S \lra S/R \lra 0\:;
  \end{equation}
  that is, $S/R$ is a flat $R$-module.

  Let $I$ be an injective $R$-module. The induced sequences
  \begin{gather*}
    0 \lra I \lra \tp{S}{I} \lra \tp{S/R}{I} \lra 0  \quad\text{and}\\
    0 \lra \Hom{S/R}{I} \lra \Hom{S}{I} \lra I \lra 0
  \end{gather*}
  are split exact so, as an $R$-module, $I$ is a direct summand of
  $\tp{S}{I}$ and of $\Hom{S}{I}$.
\end{ipg}

\begin{lem}
  \label{lem:2}
  Let $S$ be a faithfully flat $R$-algebra and let $N$ be an
  $R$-module. If $\Ext{i}{S}{N}=0$ holds for all $i>0$, then the
  following conditions are equivalent.
  \begin{eqc}
  \item $\Ext{i}{I}{N}=0$ holds for all $i>0$ and every injective
    $R$-module $I$.
  \item $\Ext[S]{i}{J}{\Hom{S}{N}}=0$ holds for all $i>0$ and every
    injective $S$-module $J$.
  \end{eqc}
\end{lem}

\begin{prf*}
  Let $X$ be an $S$-module; for every $i>0$ there are isomorphisms
  \begin{align*}
    \Ext[S]{i}{X}{\Hom{S}{N}}
    & \deq \HH[i]{\RHom[S]{X}{\Hom{S}{N}}}\\
    & \dis \HH[i]{\RHom[S]{X}{\RHom{S}{N}}}\\
    & \dis \HH[i]{\RHom{\Ltp[S]{S}{X}}{N}}\\
    & \dis \HH[i]{\RHom{X}{N}}\\
    & \deq \Ext{i}{X}{N}\:;
  \end{align*}
  the first isomorphism follows by the vanishing of $\Ext{>0}{S}{N}$,
  and the second is Hom-tensor adjointness in the derived category.
  Under the assumptions, every injective $S$-module is injective as an
  $R$-module, so it is evident from the computation that \eqclbl{i}
  implies \eqclbl{ii}.  For the converse, let $I$ be an injective
  $R$-module and recall from \pgref{ff} that it is a direct summand of
  the injective $S$-module $\Hom{S}{I}$.
\end{prf*}

For later application, we recall a fact about cotorsion modules.

\begin{ipg}
  \label{cotorsion}
  Let $X$ be a cotorsion $R$-module. For every flat $R$-module $F$ it
  follows by Hom-tensor adjointness that $\Hom{F}{X}$ is cotorsion.
\end{ipg}

We recall the notion of an injective precover, also known as an
injective right approximation.

\begin{ipg}
  Let $M$ be an $R$-module. A homomorphism $\mapdef{\f}{E}{M}$ is an
  injective \emph{precover} of $M$, if $E$ is an injective $R$-module
  and every homomorphism from an injective $R$-module to $M$ factors
  through $\f$.  Every $R$-module has an injective precover if and
  only if $R$ is noetherian; see \thmcite[2.5.1]{rha}.
\end{ipg}

In the proof of the next theorem, the noetherian hypothesis on $R$ is
used to ensure the existence of injective precovers.

\begin{thm}
  \label{thm:2}
  Let $R$ be commutative noetherian and let $S$ be a faithfully flat
  $R$-algebra. An $R$-module $N$ is Gorenstein injective if and only
  if it is cotorsion and the $S$-module $\Hom{S}{N}$ is Gorenstein
  injective.
\end{thm}

\begin{prf*}
  The ``only if'' part of the statement follows from \pgref{Jst} and
  \lemref{giascent}.  For the converse, note that \lemref{2} yields
  \begin{equation}
    \label{eq:abc}
    \Ext{>0}{I}{N}=0 \text{ holds for every injective $R$-module $I$.}
  \end{equation}
  To prove that $N$ is G-injective, it is now sufficient to show that
  $N$ has a proper left injective resolution; see \pgref{g-injective}.

  As $N$ is cotorsion, application of $\Hom{-}{N}$ to \eqref{seq}
  yields an exact sequence $0 \to \Hom{S/R}{N} \to \Hom{S}{N} \to
  \Hom{R}{N} \to 0$.  There is also an exact sequence of $S$-modules
  $0 \to G \to U \to \Hom{S}{N} \to 0$ where $U$ is injective and $G$
  is G-injective.  The composite $U \to \Hom{S}{N} \to N$ is
  surjective. Since $R$ is noetherian, there exists an injective
  precover $\mapdef{\f}{E}{N}$. The module $U$ is injective over $R$,
  so a surjective homomorphism factors through $\f$, whence $\f$ is
  surjective. Now consider the commutative diagram
  \begin{equation}
    \label{eq:diag5}
    \begin{gathered}
      \xymatrix@=1.5pc{
        &&& 0 \ar[d] \\
        &&& \Hom{S/R}{N} \ar[d] \\
        0 \ar[r] & G \ar[r] \ar[d] & U \ar[r] \ar[d]
        & \Hom{S}{N} \ar[r] \ar[d] & 0 \\
        0 \ar[r] & K \ar[r] & E \ar[r]^-{\f} & N \ar[d] \ar[r] & 0 \\
        &&& 0 }
    \end{gathered}
  \end{equation}
  To construct a proper left injective resolution of $N$, it suffices
  to show that $K$ has the same properties as $N$, including
  \eqref{abc}; that is, $K$ is cotorsion, $\Hom{S}{K}$ is G-injective
  over $S$, and $\Ext{>0}{I}{K}=0$ holds for every injective
  $R$-module $I$.

  Injective $S$-modules are injective over $R$, so by \pgref{Jst} the
  G-injective $S$-module $G$ is cotorsion over $R$. Further,
  $\Hom{S/R}{N}$ is cotorsion as $S/R$ is flat; see \pgref{cotorsion}.
  Let $F$ be a flat $R$-module; application of $\Hom{F}{-}$ to
  \eqref{diag5} yields a commutative~diagram
  \begin{equation*}
    \begin{gathered}
      \xymatrix@=1.5pc{ \Hom{F}{\Hom{S}{N}} \ar[d] \ar[r]
        & 0 \ar[d] \\
        \Hom{F}{N} \ar[d] \ar[r] & \Ext{1}{F}{K} \ar[r]
        & 0 \\
        0 \xycomma[,] }
    \end{gathered}
  \end{equation*}
  from which one concludes $\Ext{1}{F}{K}=0$. That is, $K$ is
  cotorsion.

  Let $I$ be an injective $R$-module. One has $\Ext{>0}{I}{E}=0$ and
  $\Ext{>0}{I}{N}=0$ by \eqref{abc}, so the exact sequence in
  cohomology associated to $0 \to K \to E \to N \to 0$, yields
  $\Ext{>1}{I}{K}=0$. Finally, $\Hom{I}{-}$ leaves the sequence exact,
  because $E\to N$ is an injective precover, and $\Ext{1}{I}{K}=0$
  follows.

  The exact sequence $0 \to \Hom{S}{K} \to \Hom{S}{E} \to \Hom{S}{N}
  \to 0$ shows that $\Hom{S}{K}$ has finite Gorenstein injective
  dimension over $S$, at most~$1$. \lemref{2} yields
  $\Ext{>0}{J}{\Hom{S}{K}}=0$ for every injective $S$-module $J$, so
  $\Hom{S}{K}$ is G-injective by \thmcite[2.22]{HHl04a}.
\end{prf*}

\begin{lem}
  \label{lem:1}
  Let $R$ be commutative noetherian. Let $0 \to L \to Q \to C \to 0$
  be a pure exact sequence of $R$-modules with $L\ne 0$ free and $Q$
  faithfully flat. For an $R$-module $N$ with $\Ext{i}{Q}{N}=0$ for
  all $i>0$ the next conditions are~equivalent.
  \begin{eqc}
  \item $\Ext{i}{I}{N}=0$ holds for all $i>0$ and every injective
    $R$-module $I$.
  \item $\Ext{i}{I}{\Hom{Q}{N}}=0$ holds for all $i > 0$ and every
    injective $R$-module $I$.
  \end{eqc}
  Moreover, if $N$ is cotorsion, then these conditions are equivalent
  to
  \begin{eqc}\setcounter{eqc}{2}
  \item $\Ext{i}{I}{\Hom{F}{N}}=0$ holds for all $i>0$, every
    injective $R$-module $I$, and every flat $R$-module~$F$.
  \end{eqc}
\end{lem}

\begin{prf*}
  Let $F$ be a flat $R$-module with $\Ext{>0}{F}{N}=0$. For each
  integer $i>0$ and $R$-module $X$ one has
  \begin{align*}
    \Ext{i}{X}{\Hom{F}{N}}
    & \deq \HH[i]{\RHom{X}{\Hom{F}{N}}}\\
    & \dis \HH[i]{\RHom{X}{\RHom{F}{N}}}\\
    & \dis \HH[i]{\RHom{\Ltp{F}{X}}{N}}\\
    & \dis \HH[i]{\RHom{\tp{F}{X}}{N}}\\
    & \deq \Ext{i}{\tp{F}{X}}{N} \;;
  \end{align*}
  the first isomorphism holds by the assumption on $F$, the second by
  Hom-tensor adjointness in the derived category, and the third by
  flatness of $F$.

  Let $I$ be an injective $R$-module. Since $R$ is noetherian, the
  module $\tp{F}{I}$ is injective, so it is immediate from the
  computation above that \eqclbl{i} implies \eqclbl{ii} and, if $N$ is
  cotorsion, also the stronger statement \eqclbl{iii}.

  To show that \eqclbl{ii} implies \eqclbl{i}, let $I$ be an injective
  $R$-module and consider the exact sequence
  $$0 \to \tp{L}{I} \to \tp{Q}{I} \to \tp{C}{I} \to 0.$$
  As $R$ is noetherian, the module $\tp{L}{I}$ is injective; thus the
  sequence splits, whence $I$ is isomorphic to a direct summand of
  $\tp{Q}{I}$. The computation above yields
  $\Ext{>0}{\tp{Q}{I}}{N}=0$, and $\Ext{>0}{I}{N}=0$ follows as Ext
  functors are additive.
\end{prf*}

The next result improves \cite[Ascent table I.(d)]{LWCHHl09}.  Part
\eqclbl{iii} applies, in particular, to the setting where $Q$ is a
faithfully flat $R$-algebra; cf.~\eqref{seq}.

\begin{prp}
  \label{prp:1}
  Let $R$ be commutative noetherian. For an $R$-module $N$ the
  following conditions are equivalent.
  \begin{eqc}
  \item $N$ is Gorenstein injective.
  \item $\Hom{F}{N}$ is Gorenstein injective for every flat $R$-module
    $F$.
  \item $N$ is cotorsion and $\Hom{Q}{N}$ is Gorenstein injective for
    some faithfully flat $R$-module $Q$ that contains a non-zero free
    $R$-module as a pure submodule.
  \end{eqc}
\end{prp}

\begin{prf*}
  Assume first that $N$ is G-injective and let $U$ be a totally
  acyclic complex of injective modules with $\HCy[0]{U}\is N$. Let $F$
  be a flat module. G-injective modules are cotorsion, see
  \pgref{Jst}, so the complex $\Hom{F}{U}$ is acyclic. Moreover, it is
  a complex of injective modules, and for every injective module $I$
  the complex
  \begin{equation*}
    \Hom{I}{\Hom{F}{U}} \dis \Hom{\tp{F}{I}}{U}
  \end{equation*}
  is acyclic, as $\tp{F}{I}$ is injective by the assumption that $R$
  is noetherian. This proves the implication \proofofimp[]{i}{ii}.

  That \eqclbl{ii} implies \eqclbl{iii} is trivial, just take $F=R=Q$
  and recall \pgref{Jst}. To prove that \eqclbl{iii} implies
  \eqclbl{i}, let $Q$ be a faithfully flat $R$-module with a free pure
  submodule $L\ne 0$, and assume that $\Hom{Q}{N}$ is G-injective.  In
  particular, $\Ext{>0}{I}{\Hom{Q}{N}}=0$ holds for every injective
  $R$-module $I$, so by \lemref{1}
  \begin{equation}
    \label{eq:in}
    \Ext{>0}{I}{N}=0 \text{ holds for every injective $R$-module $I$.}
  \end{equation}
  To prove that $N$ is G-injective, it is now sufficient to show that
  it has a proper left injective resolution; see
  \pgref{g-injective}. To this end, consider the pure exact sequence
  \begin{equation}
    \label{eq:purex}
    0 \lra L \lra Q \lra C \lra 0\:,
  \end{equation}
  and notice that the module $C$ is flat as $Q$ is flat.  Since
  $\Hom{Q}{N}$ is G-injective, there is an exact sequence of
  $R$-modules $0 \to G \to E \to \Hom{Q}{N} \to 0$, where $E$ is
  injective and $G$ is G-injective.  As $N$ is cotorsion, applying
  $\Hom{-}{N}$ to \eqref{purex} yields another exact sequence, and the
  two meet in the commutative diagram,
  \begin{equation}
    \label{eq:diag1}
    \begin{gathered}
      \xymatrix@=1.5pc{
        & & & 0 \ar[d] \\
        & & & \Hom{C}{N} \ar[d] \\
        0 \ar[r] & G \ar[d] \ar[r] & E \ar@{=}[d] \ar[r]
        & \Hom{Q}{N} \ar[d] \ar[r] & 0\\
        0 \ar[r] & H \ar[r] & E \ar[r]
        & \Hom{L}{N} \ar[d] \ar[r] & 0\\
        & & & 0 \xycomma[.] }
    \end{gathered}
  \end{equation}
  The first step towards construction of a proper left injective
  resolution of $N$ is to show that $H$ has the same properties as
  $N$, including \eqref{in}.

  \vspace{1ex} \emph{Claim 1:} The module $H$ is cotorsion,
  $\Hom{Q}{H}$ is G-injective, and one has $\Ext{>0}{I}{H}=0$ for
  every injective $R$-module $I$.

  \emph{Proof:} Since the module $G$ is G-injective, it is also
  cotorsion, and so is the module $\Hom{C}{N}$ as $C$ is flat; see
  \pgref{Jst} and \pgref{cotorsion}. Let $F$ be a flat $R$-module and
  apply $\Hom{F}{-}$ to \eqref{diag1} to get the commutative diagram
  \begin{equation*}
    \begin{gathered}
      \xymatrix@=1.5pc{ \Hom{F}{\Hom{Q}{N}} \ar[d] \ar[r]
        & 0 \ar[d] \\
        \Hom{F}{\Hom{L}{N}} \ar[d] \ar[r] & \Ext{1}{F}{H} \ar[r]
        & 0 \\
        0 \xycomma[,] }
    \end{gathered}
  \end{equation*}
  which shows that $\Ext{1}{F}{H}$ vanishes, whence $H$ is cotorsion.

  Notice that one has $\Hom{L}{N} \is N^\Lambda$ for some index set
  $\Lambda$. Let $I$ be an injective module; by \eqref{in} one has
  $\Ext{>0}{I}{N^\Lambda}=0$, so the exact sequence in cohomology
  induced by the last non-zero row of \eqref{diag1} yields
  $\Ext{>1}{I}{H}=0$. Now we argue that also $\Ext{1}{I}{H}$
  vanishes. As $G$ is G-injective, one has $\Ext{>0}{I}{G}=0$; by the
  assumptions on $N$ and \eqref{in}, \lemref{1} applies to yield
  $\Ext{>0}{I}{\Hom{C}{N}}=0$. Applying $\Hom{I}{-}$ to \eqref{diag1}
  one thus gets the commutative diagram
  \begin{equation}
    \label{eq:diag3}
    \begin{gathered}
      \xymatrix@=1.5pc{ \Hom{I}{\Hom{Q}{N}} \ar[d] \ar[r] &
        0 \ar[d] \\
        \Hom{I}{\Hom{L}{N}} \ar[d] \ar[r] & \Ext{1}{I}{H} \ar[r]
        & 0 \\
        0 \xycomma[,] }
    \end{gathered}
  \end{equation}
  which yields $\Ext{1}{I}{H}=0$.

  It now follows from \lemref{1} that $\Ext{>0}{I}{\Hom{Q}{H}}=0$
  holds for every injective $R$-module $I$. To prove that the module
  $\Hom{Q}{H}$ is G-injective it is thus, by \thmcite[2.22]{HHl04a},
  enough to show that it has finite Gorenstein injective
  dimension. However, that is immediate as application of $\Hom{Q}{-}$
  to last non-zero row of \eqref{diag1} yields the exact sequence
  \begin{equation*}
    0 \lra \Hom{Q}{H} \lra \Hom{Q}{E} \lra \Hom{Q}{N^\Lambda} \lra 0\;,
  \end{equation*}
  where $\Hom{Q}{E}$ is injective, and $\Hom{Q}{N^\Lambda} \is
  \Hom{Q}{N}^\Lambda$ is G-injective; see \thmcite[2.6]{HHl04a}. This
  finishes the proof of Claim 1.

  \vspace{1ex} In the commutative diagram below, the second non-zero
  row is the last non-zero row from \eqref{diag1}, and $\pi$ is a
  canonical projection.
  \begin{equation}
    \label{eq:diag4}
    \begin{gathered}
      \xymatrix@=1.5pc{
        & & & 0 \ar[d] \\
        & & & N^{\Lambda'} \ar[d] \\
        0 \ar[r] & H \ar[d] \ar[r] & E \ar[r] \ar@{=}[d]
        & N^\Lambda \ar[d]^-\pi \ar[r] & 0\\
        0 \ar[r] & K \ar[r] & E \ar[r]
        & N \ar[d] \ar[r] & 0\\
        & & & 0 }
    \end{gathered}
  \end{equation}
  To construct a proper left injective resolution of $N$, it is
  sufficient to prove that $K$ has the same properties as $N$,
  including \eqref{in}.

  \vspace{1ex}%
  \emph{Claim 2:} The module $K$ is cotorsion, $\Hom{Q}{K}$ is
  G-injective, and one has $\Ext{>0}{I}{K}=0$ for every injective
  $R$-module $I$.

  \emph{Proof:} Apply the Snake Lemma to the diagram \eqref{diag4} to
  get the exact sequence $0 \to H \to K \to N^{\Lambda'} \to 0$. By
  assumption $N$ is cotorsion, and $H$ is cotorsion by Claim 1; it
  follows that $K$ is cotorsion. Moreover, by \eqref{in} and Claim 1
  one has $\Ext{>0}{I}{H} = 0 = \Ext{>0}{I}{N}$ and hence
  $\Ext{>0}{I}{K}=0$ for every injective module~$I$.  The induced
  sequence
  \begin{equation*}
    0 \lra \Hom{Q}{H} \lra \Hom{Q}{K} \lra \Hom{Q}{N}^{\Lambda'} \lra 0
  \end{equation*}
  is exact, as $H$ is cotorsion. The modules $\Hom{Q}{N}^{\Lambda'}$
  and $\Hom{Q}{H}$ are G-injective by assumption and Claim 1, so by
  \thmcite[2.6]{HHl04a} also $\Hom{Q}{K}$ is G-injective. That
  finishes the proof of Claim 2.
\end{prf*}

\begin{cor}
  \label{cor:gidascent}
  Let $R$ be commutative noetherian. For every flat $R$-module $F$ and
  for every $R$-complex $N$ there is an inequality $\Gid{\RHom{F}{N}}
  \le \Gid{N}$.
\end{cor}

\begin{prf*}
  We may assume that $\H{N}$ is non-zero, otherwise there is nothing
  to prove.  Let $N \qra I$ be a semi-injective resolution and assume
  that $\Gid{N}=g$ is finite. By \prpref{ovl} the module $Z =
  \HCy[g]{I}$ is G-injective, and arguing as in the proof of
  \prpref[]{gidascent} one gets $\HH[i]{\RHom{F}{N}} =0$ for all
  $i>g$. The complex, $\Hom{F}{I}$ is semi-injective and isomorphic to
  $\RHom{F}{N}$ in the derived category. From \prpref{1} it follows
  that the module $\HCy[g]{\Hom{F}{I}} \is \Hom{F}{\HCy[g]{I}}$ is
  G-injective, and then \prpref{ovl} yields $\Gid{\RHom{F}{N}} \le g$.
\end{prf*}

\begin{ipg}
  \label{splf}
  For the statement of the next corollary we recall the invariant
  \begin{equation*}
    \splfR \deq \sup\setof{\pd{F}}{F\ \text{\rm is a flat $R$-module}}\:.
  \end{equation*}
  If $R$ is noetherian of finite Krull dimension $d$, then one has
  $\splfR \le d$ by works of Jensen \prpcite[6]{CUJ70}, and Gruson and
  Raynaud \thmcite[II.(3.2.6)]{LGrMRn71}. If $R$ has cardinality at
  most $\aleph_n$ for some integer $n$, then a theorem of Gruson and
  Jensen \thmcite[7.10]{LGrCUJ81} yields $\splfR \le n+1$. Osofsky
  \cite[3.1]{BLO70} gives examples of rings for which the splf
  invariant is infinite.

  Since an arbitrary direct sum of flat $R$-modules is flat, the
  invariant $\splfR$ is finite if and only if every flat $R$-module
  has finite projective dimension.
\end{ipg}

\begin{prp}
  \label{prp:2}
  Let $R$ be commutative noetherian with $\splfR < \infty$. Let $S$ be
  a faithfully flat $R$-algebra and $N$ be an $R$-module with
  $\Ext{i}{S}{N}=0$ for all $i > 0$. If\, $\Hom{S}{N}$ is Gorenstein
  injective as an $R$-module or as an $S$-module, then $N$ is
  Gorenstein injective.
\end{prp}
\enlargethispage*{\baselineskip}
\begin{prf*}
  Set $n = \splfR$.  Let $N \qra E$ be an injective resolution. As
  every flat $R$-module has projective dimension at most $n$, it
  follows by dimension shifting that the cosyzygy $Z = \HCy[n]{E}$ is
  cotorsion. Vanishing of $\Ext{>0}{S}{N}$ implies that $\Hom{S}{-}$
  leaves the sequence $0 \to N \to E^0 \to \cdots \to E^{n-1} \to Z
  \to 0$ exact.

  The class of G-injective modules is injectively resolving; see
  \thmcite[2.6]{HHl04a}. If $\Hom{S}{N}$ is G-injective over $R$, it
  thus follows that $\Hom{S}{Z}$ is G-injective over $R$. From
  \prpref{1} it follows that $Z$ is G-injective over $R$, so $N$ has
  finite Gorenstein injective dimension. Finally, \lemref{1} yields
  $\Ext{>0}{I}{N}=0$ for every injective $R$-module $I$, so $N$ is
  G-injective by \thmcite[2.22]{HHl04a}.

  If $\Hom{S}{N}$ is G-injective over $S$, then it follows as above
  that $\Hom{S}{Z}$ is G-injective over $S$. From \thmref{2} it
  follows that $Z$ is G-injective over $R$, so $N$ has finite
  Gorenstein injective dimension. Finally, it follows from \lemref{2}
  that $\Ext{>0}{I}{N}=0$ holds for every injective $R$-module $I$, so
  $N$ is G-injective.
\end{prf*}

\begin{lem}
  \label{lem:gid}
  Let $R$ be commutative noetherian and let $S$ be a faithfully flat
  $R$-algebra. For every $R$-complex $N$ of finite Gorenstein
  injective dimension one has
  \begin{equation*}
    \Gid[S]{\RHom{S}{N}} \deq \Gid{N} \deq \Gid{\RHom{S}{N}}\;.
  \end{equation*}
\end{lem}

\begin{prf*}
  We may assume that $\H{N}$ is non-zero, otherwise there is nothing
  to prove. Set $g = \Gid{N}$; \prpref{gidascent} and
  \corref{gidascent} yield $\Gid[S]{\RHom{S}{N}} \le g$ and
  $\Gid{\RHom{S}{N}} \le g$. Let $N\qra I$ be a semi-injective
  resolution; the complex $\Hom{S}{I}$ is isomorphic to $\RHom{S}{N}$
  in the derived category, and it is semi-injective as an $R$-complex
  and as an $S$-complex. By \thmcite[2.4]{JAsSSl06a} there exists an
  injective $R$-module $E$ with $\Ext{g}{E}{N} \ne 0$.

  As an $R$-module, $E$ is a direct summand of the injective
  $S$-module $\Hom{S}{E}$, see \pgref{ff}, so one gets
  \begin{align*}
    \Ext[S]{g}{\Hom{S}{E}}{\RHom{S}{N}} &\is \HH[g]{\Hom[S]{\Hom{S}{E}}{\Hom{S}{I}}}\\
    &\is \HH[g]{\Hom{\Hom{S}{E}}{I}}\\
    &\is \Ext{g}{\Hom{S}{E}}{N}\ne 0\:,
  \end{align*}
  and hence $\Gid[S]{\RHom{S}{N}} \ge g$ by \thmcite[2.4]{JAsSSl06a};
  i.e.\ equality holds.

  As an $R$-module, $E$ is a direct summand of $\tp{S}{E}$, see
  \pgref{ff}, so one gets
  \begin{align*}
    \Ext{g}{E}{\RHom{S}{N}} &\is \HH[g]{\Hom{E}{\Hom{S}{I}}}\\
    &\is \HH[g]{\Hom{\tp{S}{E}}{I}}\\
    &\is \Ext{g}{\tp{S}{E}}{N}\ne 0\:,
  \end{align*}
  and hence $\Gid{\RHom{S}{N}} \ge g$ by \thmcite[2.4]{JAsSSl06a};
  i.e.\ equality holds.
\end{prf*}

\begin{thm}
  \label{thm:3}
  Let $R$ be commutative noetherian with $\splfR < \infty$ and let $S$
  be a faithfully flat $R$-algebra. For every $R$-complex $N$ with
  $\HH[i]{N}=0$ for all $i\gg 0$ there are equalities
  \begin{equation*}
    \Gid[S]{\RHom{S}{N}} \deq \Gid{N} \deq \Gid{\RHom{S}{N}}\;.
  \end{equation*}
\end{thm}

\begin{prf*}
  In view of \lemref{gid} it is sufficient to prove that $\Gid{N}$ is
  finite if $\RHom{S}{N}$ has finite Gorenstein injective dimension as
  an $R$-complex or as an $S$-complex.  Assume, without loss of
  generality, that $\HH[i]{N}=0$ holds for $i > 0$, and let $N \qra I$
  be a semi-injective resolution.  The complex $\Hom{S}{I}$ is
  isomorphic to $\RHom{S}{N}$ in the derived category, and it is
  semi-injective as an $R$-complex and as an $S$-complex.

  If $\Gid[S]{\RHom{S}{N}}$ is finite, then one has
  $\Gid[S]{\RHom{S}{N}} \le g$ for some $g>0$. By \prpref{ovl} the
  $S$-module $\HCy[g]{\Hom{S}{I}} \is \Hom{S}{\HCy[g]{I}}$ is
  G-injective and one has $\Ext{n}{S}{\HCy[g]{I}} =
  \HH[g+n]{\Hom{S}{I}} = 0$ for all $n>0$. Now follows from \prpref{2}
  that $\HCy[g]{I}$ is G-injective, so $\Gid{N} \le g$ holds.

  The same argument applies if $\Gid{\RHom{S}{N}}$ is finite.
\end{prf*}

The equivalent of this theorem in absolute homological algebra is
\thmcite[2.2]{LWCFKk16} about injective dimension of a complex $N$. In
that statement, the \emph{a priori} assumption of vanishing of
$\HH[i]{N}$ for $i\gg 0$ is missing, though it is applied in the proof
in the same manner as above. However, \thmcite[2.2]{LWCFKk16} is
correct as stated, and an even stronger result is proved in
\corcite[3.1]{LWCSBI}.

\begin{rmk}
  \label{rmk:1}
  \thmref{2} improves \lemcite[1.3.(a)]{LWCSSW10} by removing
  assumptions that $S$ should be commutative and noetherian with
  $\dim{S}<\infty$.

  The first equality in \thmref{3} compares to the equality in
  \thmcite[1.7]{LWCSSW10} as follows: There is no assumption of finite
  Krull dimension of $R$, and $S$ is not assumed to be commutative nor
  noetherian. The inequality in \thmcite[1.7]{LWCSSW10} is subsumed by
  \corcite[9]{LZhRWe14} and compares to \prpref{gidascent} as
  discussed there.

  In \thmcite[11]{LZhRWe14} the first equality in \thmref{3} is proved
  without conditions on the homology of $N$ but under the assumption
  that $S$ is finitely generated as an $R$-module. Note that this
  assumption implies that $R$ is an algebra retract of $S$.
\end{rmk}

\section{Gorenstein flat dimension}

\noindent Recall the definition \dfncite[10.3.1]{rha} of a Gorenstein
flat module.

\begin{ipg}
  A complex $T$ of flat $R$-modules is called \emph{totally acyclic}
  if it is acyclic and the complex $\tp{I}{T}$ is acyclic for every
  injective $\Rop$-module $I$.

  An $R$-module $G$ is called \emph{Gorenstein flat} (for short,
  G-flat) if there exists a totally acyclic complex $T$ of flat
  $R$-modules with $\Co[0]{T} \is G$.
\end{ipg}

\begin{ipg}
  \label{g-flat}
  Every module has a projective resolution, so to prove that an
  $R$-module $M$ is Gorenstein flat it suffices to verify the
  following:
  \begin{rqm}
  \item $\Tor{i}{I}{M}=0$ holds for all $i>0$ and every injective
    $\Rop$-module $I$.
  \item There exists an acyclic complex $T^+ = 0 \to M \to T_{-1} \to
    T_{-2} \to \cdots$ of $R$-modules with each $T_i$ flat, such that
    $\tp{I}{T^+}$ is acyclic for every injective $\Rop$-module $I$.
  \end{rqm}
\end{ipg}

The next result is a non-commutative version of \cite[Ascent table
I.(a)]{LWCHHl09}.

\begin{lem}
  \label{lem:gfascent}
  Let $R\to S$ be a ring homomorphism such that $\fd{S}$ and
  $\fd[\Rop\!]{S}$ are finite. For a Gorenstein flat $R$-module $G$,
  the $S$-module $\tp{S}{G}$ is Gorenstein flat.
\end{lem}

\begin{prf*}
  Let $T$ be a totally acyclic complex of flat $R$-modules with
  $\Co[0]{T} \is G$. The complex $\tp{S}{T}$ consists of flat
  $S$-modules, and it is acyclic as $\fd[\Rop\!]{S}$ is finite; see
  \lemcite[2.3]{CFH-06}. One has $\Co[0]{\tp{S}{T}} \is \tp{S}{G}$, so
  it suffices to show that the complex $\tp[S]{J}{\tpp{S}{T}} \is
  \tp{J}{T}$ is acyclic for every injective $\Sop$-module $J$. As
  $\fd{S}$ is finite, every injective $\Sop$-module has finite
  injective dimension over $\Rop$; see
  \corcite[4.2]{LLAHBF91}. Acyclicity of $\tp{J}{T}$ hence follows
  from \lemcite[2.3]{CFH-06}.
\end{prf*}

In standard homological algebra, the study of flat modules and flat
dimension is propelled by flat--injective duality: An $R$-module $F$
is flat if and only if the $\Rop$-module $\Hom[\ZZ]{F}{\QQ/\ZZ}$ is
injective. A similar Gorenstein flat--injective duality is only known
hold if $R$ is right coherent, so that is a convenient setting for
studies of Gorenstein flat dimension; see for example
\thmcite[3.14]{HHl04a}.

As Bennis demonstrates \thmcite[2.8]{DBn09}, one can do a little
better and get a well-behaved theory of Gorenstein flat dimension, as
long as the class of Gorenstein flat modules is projectively
resolving. Rings with that property have become known as left
GF-closed; we recall the definition in \pgref{GF-closed}.

Iacob \cite{AIc09} defined a notion of Gorenstein flat dimension for
unbounded complexes, but only for complexes over left GF-closed
rings. Here we give a definition, inspired by
Veliche's~\dfncite[3.1]{OVl06}; it applies to complexes over any ring
and coincides with Iacob's for left GF-closed rings.

To get started we recall some terminology from Liang~\cite{LLn13}.

\begin{ipg}
  Let $M$ be an $R$-complex; a \emph{semi-flat replacement} of $M$ is
  a semi-flat $R$-complex $F$ such that there is an isomorphism $F
  \qis M$ in the derived category.  Every complex has a
  semi-projective resolution and hence a semi-flat replacement; see
  \pgref{res}.  A \emph{Tate flat resolution} of $M$ is a pair $(T,F)$
  where $T$ a totally acyclic complex of flat $R$-modules and $F\simeq
  M$ a semi-flat replacement with $\Thb{g}{T} \is \Thb{g}{F}$ for some
  $g\in\mathbb{Z}$.
\end{ipg}

\begin{dfn}
  \label{dfn:gfd}
  Let $M$ be an $R$-complex. The \emph{Gorenstein flat dimension} of
  $M$ is given by
  \begin{equation*}
    \Gfd{M} \deq \inf\setof{g\in\ZZ}{(T,F)\
      \text{is a Tate flat resolution of $M$ with $\Thb{g}{T}
        \is \Thb{g}{F}$}}\:.
  \end{equation*}
\end{dfn}

\begin{ipg}
  \label{df}
  Let $M$ be an $R$-complex; the following facts are evident from the
  definition.
  \begin{prt}
  \item One has $\Gfd{M} < \infty$ if and only if $M$ admits a Tate
    flat resolution.
  \item $M$ is acyclic if and only if one has $\Gfd{M} = -\infty$.
  \item $\Gfd{(\Shift[n]{M})} = \Gfd{M}+n$ holds for every
    $n\in\ZZ$.\footnote{ Here $\Shift[n]{M}$ is the complex with
      $(\Shift[n]{M})_i = M_{i-n}$ and $\dif[i]{\Shift[n]{M}} =
      (-1)^n\dif[i-n]{M}$ for all $i\in\ZZ$.}
  \end{prt}
\end{ipg}

A Gorenstein flat $R$-module $M\ne 0$ clearly has has $\Gfd{M}=0$ per
\dfnref{gfd}. In \rmkref{iacob} we compare \dfnref[]{gfd} to competing
definitions in the literature.

\begin{prp}
  \label{prp:property}
  Let $M$ be an $R$-complexes and $g\in\ZZ$. The following conditions
  are equivalent.
  \begin{eqc}
  \item $\Gfd{M} \le g$.
  \item $\H[i]{M} =0$ holds for all $i > g$, and there exists a
    semi-flat replacement $F\qis M$ such that $\Co[i]{F}$ is
    Gorenstein flat for every $i\ge g$.
  \item $\H[i]{M} =0$ holds for all $i > g$, and there exists a
    semi-flat replacement $F\qis M$ such that $\Co[g]{F}$ is
    Gorenstein flat.
  \end{eqc}
\end{prp}

\begin{prf*}
  If $\Gfd{M}\le g$ holds, then there exists a Tate flat resolution
  $(T,F)$ of $M$ with $\Thb{g}{T} \is \Thb{g}{F}$. As $T$ is acyclic
  and $F \qis M$ one has $\H[i]{M} \is \H[i]{F} \is \H[i]{T} = 0$ for
  all $i > g$, and $\Co[i]{F} \is \Co[i]{T}$ is G-flat for every $i\ge
  g$. Thus \eqclbl{i} implies \eqclbl{ii}.

  Clearly, \eqclbl{ii} implies \eqclbl{iii}; to finish the proof
  assume that $\H[i]{M} = 0$ holds for all $i > g$ and that there
  exists a semi-flat replacement $F\qis M$ such that $\Co[g]{F}$ is
  G-flat. As one has $\H{F} \is \H{M}$ it follows that the complex
  \begin{equation*}
    F' \deq \cdots \lra F_{g+1} \lra F_{g} \lra \Co[g]{F} \lra 0
  \end{equation*}
  is acyclic; i.e.\ it is an augmented flat resolution of
  $\Co[g]{F}$. For every injective $\Rop$-module $I$ one has
  $\Tor{>0}{I}{\Co[g]{F}}=0$, so the complex $\tp{I}{F'}$ is
  acyclic. Moreover, there exists an acyclic complex of $R$-modules
  \begin{equation*}
    T' \deq 0 \lra \Co[g]{F} \lra T_{g-1} \lra T_{g-2} \lra \cdots
  \end{equation*}
  with each $T_{i}$ flat, such that $\tp{I}{T'}$ is acyclic for every
  injective $\Rop$-module $I$.  Now let $T$ be the complex obtained by
  splicing together $F'$ and $T'$ at $\Co[g]{F}$; that is
  \begin{equation*}
    T \deq \cdots \lra F_{g+1} \lra F_{g} \lra T_{g-1} \lra T_{g-2} \lra \cdots\: .
  \end{equation*}
  It is a totally acyclic complex of flat $R$-modules, so $(T,F)$ is a
  Tate flat resolution of $M$ with $\Thb{g}{T} \is \Thb{g}{F}$, and so
  one has $\Gfd{M} \le g$.
\end{prf*}

In the next statement `fd' means flat dimension as defined in
\seccite[2.F]{LLAHBF91}.

\begin{prp}
  \label{prp:fd}
  Let $M$ be an $R$-complex; the following assertions hold.
  \begin{prt}
  \item One has $\Gfd{M} \le \fd{M}$; equality holds if
    $\fd{M}<\infty$.
  \item One has $\Gid[\Rop\!]{\Hom[\ZZ]{M}{\QQ/\ZZ}} \le \Gfd{M}$;
    equality holds if $R$ is~ right coherent.
  \end{prt}
\end{prp}

\begin{prf*}
  $(a)$: For an acyclic $R$-complex $M$ one has $\Gfd{M} = -\infty =
  \fd{M}$, and the inequality is trivial if $\fd{M} = \infty$. Now
  assume that $\H{M} \ne 0$ and $\fd{M}=g$ holds for some
  $g\in\ZZ$. It follows from \thmcite[2.4.F]{LLAHBF91} that $\H[i]{M}
  = 0$ holds for all $i > g$, and that there exists a semi-flat
  replacement $F\qis M$ such that $\Co[g]{F}$ is flat, in particular
  G-flat. Thus one has $\Gfd{M} \le g$ by \prpref{property}. To prove
  that equality holds, assume towards a contradiction that $\Gfd{M}\le
  g-1$ holds. One then has $\H[i]{M} =0$ for all $i > g-1$ and there
  exists a semi-flat replacement $F\qis M$ such that $\Co[g-1]{F}$ is
  G-flat; see \prpref{property}. Consider the exact sequence
  \begin{equation*}
    0 \lra \Co[g]{F} \lra F_{g-1} \lra \Co[g-1]{F} \lra 0\: .
  \end{equation*}
  The module $\Co[g]{F}$ is flat, so one has $\fd(\Co[g-1]{F})\leq 1$,
  and then $\Co[g-1]{F}$ is flat; see
  Bennis~\thmcite[2.2]{DBn11}. Thus one has $\fd{M}\leq g-1$, which is
  a contradiction.

  $(b)$: For an acyclic $R$-complex $M$ one has
  $\Gid[\Rop\!]{\Hom[\ZZ]{M}{\QQ/\ZZ}} = -\infty = \Gfd{M}$, and the
  inequality is trivial if $\Gfd{M} = \infty$. Now assume that
  $\H{M}\ne 0$ and $\Gfd{M} = g$ holds for some $g\in\ZZ$. By
  \prpref{property} one has $\H[i]{M} = 0$ for all $i>g$, and there
  exists a semi-flat replacement $F\qis M$ such that $\Co[g]{F}$ is
  G-flat. The $\Rop$-complex $\Hom[\ZZ]{F}{\QQ/\ZZ}$ is semi-injective
  and yields a semi-injective resolution $\Hom[\ZZ]{M}{\QQ/\ZZ} \qra
  \Hom[\ZZ]{F}{\QQ/\ZZ}$; see~\cite[1.4.I]{LLAHBF91}. For all $i > g$
  one has
  \begin{equation*}
    \HH[i]{\Hom[\ZZ]{M}{\QQ/\ZZ}} \dis \Hom[\ZZ]{\H[i]{M}}{\QQ/\ZZ} = 0\:,
  \end{equation*}
  and by \prpcite[3.11]{HHl04a} the $\Rop$-module
  $\HCy[g]{\Hom[\ZZ]{F}{\QQ/\ZZ}} \is \Hom[\ZZ]{\Co[g]F}{\QQ/\ZZ}$ is
  G-injective, so one gets $\Gid[\Rop\!]{\Hom[\ZZ]{M}{\QQ/\ZZ}} \le g$
  from \prpref{ovl}.

  Assume now that $g = \Gid[\Rop\!]{\Hom[\ZZ]{M}{\QQ/\ZZ}}$ is
  finite. Let $F \qis M$ be a semi-flat replacement; as above it
  yields a semi-injective resolution
  \begin{equation*}
    \Hom[\ZZ]{M}{\QQ/\ZZ} \qra \Hom[\ZZ]{F}{\QQ/\ZZ}\:.
  \end{equation*}
  By \prpref{ovl} one has $\Hom[\ZZ]{\H[i]{M}}{\QQ/\ZZ} \is
  \HH[i]{\Hom[\ZZ]{M}{\QQ/\ZZ}} = 0$ for all $i > g$ and the module
  $\Hom[\ZZ]{\Co[g]{F}}{\QQ/\ZZ} \is \HCy[g]{\Hom[\ZZ]{F}{\QQ/\ZZ}}$
  is G-injective. It follows that $\H[i]{M}=0$ holds for all $i > g$,
  and under the assumption that $R$ is right coherent it follows from
  \prpcite[3.11]{HHl04a} that $\Co[g]{F}$ is G-flat. Thus one has
  $\Gfd{M} \le g$ by \prpref{property}, and
  $\Gid[\Rop\!]{\Hom[\ZZ]{M}{\QQ/\ZZ}} = \Gfd{M}$ holds.
\end{prf*}

The next result is proved by standard applications of the Snake Lemma.

\begin{lem}
  \label{lem:cz}
  Let $0 \to X \to Y \to Z \to 0$ be an exact sequence of
  $R$-complexes and let $n\in\mathbb{Z}$. The following assertions
  hold.
  \begin{prt}
  \item If\, $\H[n]{Z}=0$, then the sequence $0 \to \Co[n-1]{X} \to
    \Co[n-1]{Y} \to \Co[n-1]{Z} \to 0$ is exact. Moreover if
    $\H[n]{Y}=0$, then the converse holds.
  \item If\, $\H[n]{X}=0$, then the sequence $0 \to \Cy[n+1]{X} \to
    \Cy[n+1]{Y} \to \Cy[n+1]{Z} \to 0$ is exact. Moreover if
    $\H[n]{Y}=0$, then the converse holds.\qed
  \end{prt}
\end{lem}

\begin{ipg}
  \label{GF-closed}
  Recall from \cite{DBn09} that $R$ is called left \emph{GF-closed} if
  the class of Gorenstein flat $R$-modules is closed under extensions
  (equivalently, it is projectively resolving). The class of left
  GF-closed rings strictly contains the class of right coherent rings.
\end{ipg}

The next result appears as \lemcite[3.2]{AIc09} with the additional
hypothesis that $R$ is left GF-closed; that hypothesis is, however,
not needed.

\begin{lem}
  \label{lem:iacob}
  Let $0 \to F \to G \to H \to 0$ be an exact sequence of
  $R$-modules. If $F$ is flat, $G$ is Gorenstein flat, and
  $\Hom[\ZZ]{H}{\QQ/\ZZ}$ is Gorenstein injective, then $H$ is
  Gorenstein flat.
\end{lem}

\begin{prf*}
  It is argued in the proof of \lemcite[3.2]{AIc09} that there is an
  exact sequence $0 \to H \to X \to Y \to 0$ with $X$ flat and $Y$
  G-flat, so $H$ is G-flat by \lemcite[2.4]{DBn09}.
\end{prf*}

\begin{prp}
  \label{prp:gfd-sp}
  Let $M$ be an $R$-complex and $g\in\mathbb{Z}$.  Conditions
  \eqclbl{i} and \eqclbl{ii} below are equivalent and imply
  \eqclbl{iii}. Moreover, if $R$ is left GF-closed, then all three
  conditions are equivalent.
  \begin{eqc}
  \item $\H[i]{M}=0$ holds for all $i > g$, and for every semi-flat
    replacement $F \qis M$ the cokernel $\Co[g]{F}$ is Gorenstein
    flat.
  \item $\H[i]{M}=0$ holds for all $i > g$, and there exists a
    semi-projective resolution $P \qra M$ such that the cokernel
    $\Co[g]{P}$ is Gorenstein flat.
  \item[$(iii)$] $\Gfd{M} \le g$.
  \end{eqc}
\end{prp}

\begin{prf*}
  Clearly \eqclbl{i} implies \eqclbl{ii}, and \eqclbl{ii} implies
  \eqclbl{iii} by \prpref{property}.

  \proofofimp{ii}{i} Let $P \qra M$ be a semi-projective resolution
  such that $\Co[g]{P}$ is G-flat. Let $F \simeq M$ be a semi-flat
  replacement. There is by \cite[1.4.P]{LLAHBF91} a quasi-isomorphism
  $\mapdef[\qra]{\f}{P}{F}$, and after adding to $P$ a projective
  precover of $F$, we may assume that $\f$ is surjective. One then has
  an exact sequence $0 \to K \to P \to F \to 0$. The kernel $K$ is
  semi-flat and acyclic, so the module $\Co[i]{K}$ is flat for every
  $i\in\mathbb{Z}$; see Christensen and
  Holm~\thmcite[7.3]{LWCHHl15}. The sequence
  \begin{equation}
    \label{eq:c}
    0 \lra \Co[g]{K} \lra \Co[g]{P} \lra \Co[g]{F} \lra 0
  \end{equation}
  is exact by \lemref{cz}(a). The module $\Co[g]{P}$ is G-flat, so to
  prove that $\Co[g]{F}$ is G-flat, it is by \lemref{iacob} enough to
  prove that $\Hom[\ZZ]{\Co[g]F}{\QQ/\ZZ}$ is G-injective. By
  \prpref[Propositions~]{property} and \prpref[]{fd}(b) one has
  $\Gid[\Rop\!]{\Hom[\ZZ]{M}{\QQ/\ZZ}}\le g$. The $\Rop$-complex
  $\Hom[\ZZ]{F}{\QQ/\ZZ}$ is semi-injective and yields a
  semi-injective resolution $\Hom[\ZZ]{M}{\QQ/\ZZ} \qra
  \Hom[\ZZ]{F}{\QQ/\ZZ}$; see \cite[1.4.I]{LLAHBF91}. It now follows
  from \prpref{ovl} that the module $\Hom[\ZZ]{\Co[g]{F}}{\QQ/\ZZ} \is
  \HCy[g]{\Hom[\ZZ]{F}{\QQ/\ZZ}}$ is G-injective.

  Thus, the conditions \eqclbl{i} and \eqclbl{ii} are equivalent. We
  now assume that $R$ is left GF-closed and verify the implication
  \proofofimp[]{iii}{ii}.  By \prpref{property} one has $\H[i]{M} = 0$
  for all $i > g$, and there exists a semi-flat replacement $F\qis M$
  with $\Co[g]{F}$ G-flat. Choose a surjective semi-projective
  resolution $P \qra F$; it yields by \cite[1.4.P]{LLAHBF91} a
  semi-projective resolution $P \qra M$. As above one gets an sequence
  \eqref{c}, where now $\Co[g]{K}$ is flat and $\Co[g]{F}$ is G-flat;
  it follows that $\Co[g]{P}$ is G-flat as $R$ is left GF-closed.
\end{prf*}

For commentary on the proof, see \rmkref{driss} at the end of the
section.

\begin{rmk}
  \label{rmk:iacob}
  From \prpref{gfd-sp} and \thmcite[1]{AIc09} it is clear that for
  complexes over a left GF-closed ring, \dfnref{gfd} agrees with
  Iacob's definition \dfncite[3.2]{AIc09}. In this setting, it thus
  extends the definitions in \dfncite[3.9]{HHl04a} and
  \cite[1.9]{CFH-06} of Gorenstein flat dimension for modules and
  complexes with bounded below homology; see \rmkcite[3]{AIc09}.

  Though they potentially differ over rings that are not GF-closed, we
  use the same symbol, `Gfd', for our notion and the one from
  \cite[1.9]{CFH-06}. It causes no ambiguity, as the results we prove
  in the remainder of this paper are valid for either notion. When
  needed, we make remarks to that effect.
\end{rmk}

The next result improves \thmcite[12]{LZhRWe14} by removing
assumptions that $R$ and $S$ should be commutative noetherian.

\begin{prp}
  \label{prp:gfdascent}
  Let $R\to S$ be a ring homomorphism such that $\fd{S}$ and
  $\fd[\Rop\!]{S}$ are finite. For every $R$-complex $M$ one has
  $\Gfd[S]{\Ltpp{S}{M}} \le \Gfd{M}$.
\end{prp}

\begin{prf*}
  We may assume that $\H{M}$ is non-zero, otherwise there is nothing
  to prove. Assume that $\Gfd{M} = g$ holds for some $g\in\ZZ$. By
  \prpref{property} one has $\H[i]{M}=0$ for all $i > g$, and there
  exists a semi-flat replacement $F\qis M$ with $\Co[g]{F}$
  G-flat. The $S$-complex $\tp{S}{F}$ is semi-flat and isomorphic to
  $\Ltp{S}{M}$ in the derived category. For $n>0$ one has
  \begin{equation*}
    \H[g+n]{\Ltp{S}{M}} \dis \H[g+n]{\tp{S}{F}} \dis \Tor{n}{S}{\Co[g]{F}} = 0\:,
  \end{equation*}
  where vanishing follows from \lemcite[2.3]{CFH-06} as
  $\fd[\Rop\!]{S}$ is finite. Finally, the module $\Co[g]{\tp{S}{F}}
  \is \tp{S}{\Co[g]{F}}$ is G-flat by \lemref{gfascent}, so
  $\Gfd[S]{\Ltpp{S}{M}} \le g$ holds by \prpref{property}.
\end{prf*}

\begin{rmk}
  For an $R$-complex $M$ with $\H[i]{M}=0$ for $i \ll 0$ the
  inequality in \prpref{gfdascent} is also valid with the definition
  of Gorenstein flat dimension from \cite[1.9]{CFH-06}. Indeed, the
  argument for \cite[(4.4)]{LWCHHl09} also applies when $R$ is not
  commutative.
\end{rmk}

\begin{rmk}
  \label{rmk:driss}
  From the proof of \proofofimp[]{iii}{ii} in \prpref{gfd-sp},
  consider the exact sequence
  \begin{equation}
    \label{eq:driss}
    0 \lra \Co[g]{K} \lra \Co[g]{P} \lra \Co[g]{F} \lra 0
  \end{equation}
  where the module $\Co[g]{K}$ is flat and $\Co[g]{F}$ is Gorenstein
  flat. Under the assumption that $R$ is left GF-closed, the module in
  the middle, $\Co[g]{P}$, is Gorenstein flat.  Without that
  assumption, the best we know is to consider the pull-back diagram
  \begin{equation*}
    \xymatrix@=1.5pc{
      && 0 \ar[d] & 0 \ar[d]\\
      && G \ar@{=}[r] \ar[d] & G \ar[d] \\
      0 \ar[r] & \Co[g]{K} \ar@{=}[d] \ar[r]
      & X \ar[d] \ar[r] & T \ar[d] \ar[r] & 0\\
      0 \ar[r] \ & \Co[g]{K} \ar[r] & \Co[g]{P} \ar[r] \ar[d]
      & \Co[g]{F} \ar[r] \ar[d]& 0\\
      && 0  & 0
    }
  \end{equation*}
  where $T$ is flat and $G$ is Gorenstein flat. It follows that $X$ is
  flat, so $\Co[g]{P}$ has Gorenstein flat dimension at most $1$, per
  \dfnref{gfd} as well as per \dfncite[3.9]{HHl04a}.

  Notice that the module $\Co[g]{K}$ in \eqref{driss} is not just
  Gorenstein flat but flat. In an earlier version of this paper we
  asked if a ring is GF-closed provided that in every exact sequence
  of $R$-modules $0 \to F \to M \to G \to 0$ with $F$ flat and $G$
  Gorenstein flat also $M$ is Gorenstein flat. A positive answer to
  this question is provided in the Appendix by Driss Bennis.
\end{rmk}

\section{Faithfully flat base change}

\noindent
In this section, $R$ is commutative and $S$ is an $R$-algebra.

\begin{lem}
  \label{lem:3}
  Let $S$ be a faithfully flat $R$-algebra.  For an $R$-module $M$ the
  following conditions are equivalent.
  \begin{eqc}
  \item $\Tor{i}{I}{M}=0$ holds for all $i>0$ and every injective
    $R$-module $I$.
  \item $\Tor[S]{i}{J}{\tp{S}{M}}=0$ holds for all $i>0$ and every
    injective $\Sop$-module $J$.
  \end{eqc}
\end{lem}

\begin{prf*}
  For every $\Sop$-module $X$ and every $i >0$ one has
  \begin{equation*}
    \Tor[S]{i}{X}{\tp{S}{M}} \dis \H[i]{\Ltp[S]{X}{\Ltpp{S}{M}}} \is
    \H[i]{\Ltp{X}{M}} \dis \Tor{i}{X}{M}\:.
  \end{equation*}
  As $S$ is flat over $R$, every injective $\Sop$-module is an
  injective $R$-module, so \eqclbl{i} implies \eqclbl{ii}. Let $I$ be
  an $R$-module and recall from \pgref{ff} that it is a direct summand
  of the injective $\Sop$-module $\Hom{S}{I}$; now \eqclbl{i} follows
  from \eqclbl{ii}.
\end{prf*}

For use in the next proof, we recall the notion of a flat preenvelope,
also known as a flat left approximation.

\begin{ipg}
  Let $M$ be an $R$-module. A homomorphism $\mapdef{\f}{M}{F}$ is a
  flat \emph{preenvelope} of $M$, if $F$ is a flat $R$-module and
  every homomorphism from $M$ to a flat $R$-module factors through
  $\f$.  Every $R$-module has a flat preenvelope if and only if $R$ is
  coherent; see \thmcite[5.4.1]{rha}.
\end{ipg}

\begin{lem}
  \label{lem:css}
  Let $R$ be coherent and let $S$ be a faithfully flat $R$-algebra
  that is left GF-closed. An $R$-module $M$ is Gorenstein flat if and
  only if the $S$-module $\tp{S}{M}$ is Gorenstein flat.
\end{lem}

\begin{prf*}
  The ``only if'' part is a special case of \lemref{gfascent}. To
  prove ``if'', assume that $\tp{S}{M}$ is a G-flat $S$-module. Now
  $\Tor[S]{>0}{J}{\tp{S}{M}}=0$ holds for every injective
  $\Sop$-module $J$, so by \lemref{3} one has
  \begin{equation}
    \label{eq:torvan0}
    \Tor{>0}{I}{M}=0 \text{ for every injective $R$-module $I$.}
  \end{equation}
  Thus, \pgref{g-flat}.(1) holds and we proceed to verify
  \pgref{g-flat}.(2).

  Let $\mapdef{\f}{M}{F}$ be a flat preenvelope; it exists because $R$
  is coherent. The G-flat module $\tp{S}{M}$ embeds into a flat
  $S$-module, and $M$ embeds into $\tp{S}{M}$ as $S$ is faithfully
  flat; cf.~\eqref{seq}. A flat $S$-module is also flat as an
  $R$-module, so $M$ embeds into a flat $R$-module, and it follows
  that $\f$ is injective.  Consider the associated exact sequence
  \begin{equation*}
    \eta \deq 0 \lra M \xra{\f} F \lra C \lra 0\;.
  \end{equation*}
  Let $I$ be an injective $R$-module. The induced sequence
  $\tp{I}{\eta}$ is exact if and only if
  $\Hom[\ZZ]{\tp{I}{\eta}}{\QQ/\ZZ} \is
  \Hom{\eta}{\Hom[\ZZ]{I}{\QQ/\ZZ}}$ is exact. The $R$-module
  $\Hom[\ZZ]{I}{\QQ/\ZZ}$ is flat, as $R$ is coherent, so exactness of
  $\Hom{\eta}{\Hom[\ZZ]{I}{\QQ/\ZZ}}$ holds because $\f$ is a flat
  preenvelope. Notice that exactness of $\tp{I}{\eta}$ and
  \eqref{torvan0} yield
  \begin{equation}
    \label{eq:torvan}
    \Tor{>0}{I}{C}=0 \text{ for every injective $R$-module $I$.}
  \end{equation}
  To construct the complex in \pgref{g-flat}.(2), it is now sufficient
  to show that $\tp{S}{C}$ is G-flat over $S$. Consider the following
  push-out diagram in the category of $S$-modules where $H$ is flat
  and $G$ is G-flat.
  \begin{equation*}
    \xymatrix@=1.5pc{
      & 0 \ar[d] & 0 \ar[d]\\
      0 \ar[r] & \tp{S}{M} \ar[r] \ar[d] & \tp{S}{F} \ar[r] \ar[d]
      & \tp{S}{C} \ar@{=}[d] \ar[r] & 0\\
      0 \ar[r] & H \ar[d] \ar[r] & X \ar[r] \ar[d] & \tp{S}{C} \ar[r] & 0\\
      & G \ar@{=}[r] \ar[d] & G \ar[d] \\
      & 0& 0
    }
  \end{equation*}
  Since $S$ is left GF-closed, the module $X$ is G-flat. To conclude
  from \lemref{iacob}, applied to the diagram's second non-zero row,
  that $\tp{S}{C}$ is G-flat, we need to verify that the $\Sop$-module
  $\Hom[\ZZ]{\tp{S}{C}}{\QQ/\ZZ}$ is G-injective. Application of
  $\Hom[\ZZ]{-}{\QQ/\ZZ}$ to the second row shows that the module has
  finite Gorenstein injective dimension (at most $1$). Let $J$ be an
  injective $\Sop$-module and $i > 0$ be an integer. The isomorphism
  below is Hom-tensor adjointness, and vanishing follows from
  \lemref{3} in view of \eqref{torvan}
  \begin{align*}
    \Ext[\Sop\!]{i}{J}{\Hom[\ZZ]{\tp{S}{C}}{\QQ/\ZZ}} \dis
    \Hom[\ZZ]{\Tor[S]{i}{J}{\tp{S}{C}}}{\QQ/\ZZ} = 0\:.
  \end{align*}
  It follows that $\Hom[\ZZ]{\tp{S}{C}}{\QQ/\ZZ}$ is G-injective; see
  \thmcite[2.22]{HHl04a}.
\end{prf*}

It is straightforward to verify that for a flat $R$-module $F$ and a
Gorenstein flat $R$-module $G$, the module $\tp{F}{G}$ is Gorenstein
flat; see \cite[Ascent table I.(a)]{LWCHHl09}. The next lemma provides
a partial converse; it applies, in particular, to the setting where
$Q$ is a faithfully flat $R$-algebra; cf.~\eqref{seq}.

\begin{prp}
  \label{prp:cs}
  Let $R$ be coherent and $M$ be an $R$-module. If $\tp{Q}{M}$ is
  Gorenstein flat for some faithfully flat $R$-module $Q$ that
  contains a non-zero free $R$-module as a pure submodule, then $M$ is
  Gorenstein flat.
\end{prp}

\begin{prf*}
  As $Q$ is faithfully flat, vanishing of $\Tor{i}{I}{\tp{Q}{M}}$
  implies vanishing of $\Tor{i}{I}{M}$, so one has $\Tor{>0}{I}{M}=0$
  for every injective $R$-module $I$.  Thus \pgref{g-flat}.(1) is
  satisfied and we proceed to construct the complex in
  \pgref{g-flat}.(2).

  Let $L\ne 0$ be a free pure submodule of $Q$; application of
  $\tp{-}{M}$ to the associated pure embedding yields an embedding
  $\tp{L}{M} \to \tp{Q}{M}$. By assumption, the module $\tp{Q}{M}$
  embeds into a flat $R$-module, whence $\tp{L}{M}$, and therefore
  $M$, embeds into a flat $R$-module. Since $R$ is coherent, there
  exists a flat preenvelope, $\mapdef{\f}{M}{F}$, and it is
  necessarily injective.  Consider the associated exact sequence $0
  \to M \to F \to C \to 0$. As in the proof of \lemref{css} one gets
  \begin{equation}
    \label{eq:torvan1}
    \Tor{>0}{I}{C}=0 \text{ for every injective $R$-module $I$.}
  \end{equation}
  To construct the complex in \pgref{g-flat}.(2) it is now sufficient
  to show that $\tp{Q}{C}$ is G-flat. It is immediate from the induced
  exact sequence
  \begin{equation*}
    0 \lra \tp{Q}{M} \lra \tp{Q}{F} \lra \tp{Q}{C} \lra 0
  \end{equation*}
  that $\tp{Q}{C}$ has finite Gorenstein flat dimension. For every
  injective $R$-module $I$ one gets
  $$\Tor{>0}{I}{\tp{Q}{C}} = 0$$
  from \eqref{torvan1}, so $\tp{Q}{C}$ is G-flat by
  \thmcite[3.14]{HHl04a}.
\end{prf*}

\begin{prp}
  \label{prp:gfdascent1}
  For every flat $R$-module $F$ and for every $R$-complex $M$ there is
  an inequality
  $$\Gfd{\tpp{F}{M}} \le \Gfd{M}.$$
\end{prp}

\begin{prf*}
  We may assume that $\H{M}$ is non-zero, otherwise there is nothing
  to prove. Assume that $\Gfd{M} = g$ holds for some $g\in\ZZ$. By
  \prpref{property} one has $\H[i]{M}=0$ for all $i > g$, and there
  exists a semi-flat replacement $P\qis M$ with $\Co[g]{P}$
  G-flat. The complex $\tp{F}{P}$ is semi-flat and isomorphic to
  $\tp{F}{M}$ in the derived category. For $i>g$ one has
  $\H[i]{\tp{F}{M}} \is \tp{F}{\H[i]{M}} = 0$, and the module
  $\Co[g]{\tp{F}{P}} \is \tp{F}{\Co[g]{P}}$ is G-flat by \cite[Ascent
  table I.(a)]{LWCHHl09}, so $\Gfd{\tpp{F}{M}} \le g$ holds by
  \prpref{property}.
\end{prf*}

\begin{rmk}
  For an $R$-complex $M$ with $\H[i]{M}=0$ for $i \ll 0$ the
  inequality in \prpref{gfdascent1} is by \cite[(4.4)]{LWCHHl09} valid
  for Gorenstein flat dimension defined in \cite[1.9]{CFH-06}.
\end{rmk}

\begin{thm}
  \label{thm:transfer}
  Let $R$ be commutative coherent, let $S$ be a faithfully flat
  $R$-algebra, and let $M$ be an $R$-complex. There is an equality
  \begin{equation*}
    \Gfd{M} \deq \Gfd{\tpp[R]{S}{M}}\;,
  \end{equation*}
  and if $S$ is left GF-closed also an equality
  \begin{equation*}
    \Gfd[S]{\tpp[R]{S}{M}} \deq \Gfd{M}\:.
  \end{equation*}
\end{thm}

\begin{prf*}
  We may assume that $\H{M}$ is non-zero, otherwise there is nothing
  to prove.

  Consider the first equality. The inequality ``$\ge$'' is a special
  case of \prpref{gfdascent1}.  For the opposite inequality, assume
  $\Gfd{\tpp{S}{M}} = g$ holds for some $g\in\ZZ$. Choose a
  semi-projective resolution $P \qra M$; the induced quasi-isomorphism
  $\tp{S}{P} \to \tp{S}{M}$ is a semi-flat replacement. From
  \prpref{gfd-sp} it follows that $\tp{S}{\H[i]{M}} \is
  \H[i]{\tp{S}{M}} = 0$ holds for all $i > g$ and that the module
  $\tp{S}{\Co[g]{P}} \is \Co[g]{\tp{S}{P}}$ is G-flat. Since $S$ is
  faithfully flat, it follows that $\H[i]{M} = 0$ holds for all $i>
  g$, and by \prpref{cs} the module $\Co[g]{P}$ is G-flat. By
  \prpref{gfd-sp} one now has $\Gfd{M} \le g$, and this proves the
  first equality.

  A parallel argument applies to establish the second equality; only
  one has to invoke \prpref{gfdascent} instead of
  \prpref[]{gfdascent1} and \lemref{css} instead of \prpref[]{cs}.
\end{prf*}

\begin{rmk}
  \label{rmk:2}
  The first equality in \thmref{transfer} improves the second equality
  in \prpcite[1.9]{LWCSSW10} by removing the requirements that $R$ be
  semi-local and $S$ be commutative noetherian. The equality in
  \thmcite[1.8]{LWCSSW10} are subsumed by \thmcite[15]{LZhRWe14}, and
  the second equality in \thmref{transfer} improves
  \thmcite[15]{LZhRWe14} by removing the assumption that $S$ is
  commutative and by relaxing the conditions on $R$ and $S$ from
  noetherian to coherent/left GF-closed.
\end{rmk}

\section{Closing remarks}
\noindent
Research on Gorenstein homological dimensions continues to be guided
by the idea that every result about absolute homological dimensions
should have a counterpart in Gorenstein homological algebra. In a
curious departure from this principle, the behavior of Gorenstein
injective dimension under faithfully flat base change was understood
\thmcite[1.7]{LWCSSW10} years before the absolute case
\thmcite[2.2]{LWCFKk16}. The main results of the present paper have
brought the (co-)basechange results for Gorenstein dimensions in
closer alignment with the results for absolute homological dimensions
by relaxing the conditions on rings and complexes. It would still be
interesting to remove the \emph{a priori} assumption of homological
boundedness from \thmref{3} to align it even closer with
\corcite[3.1]{LWCSBI}.

\section*{Acknowledgments}
\noindent
This paper was written during visits by L.L.\ to Texas Tech
University---while F.K.\ was a student there---and by L.W.C.\ and
L.L.\ to Nanjing University; we gratefully acknowledge the hospitality
of both institutions. In particular, we thank Nanqing Ding, our host
at Nanjing, for inviting and supporting us.

\appendix
\section*{Appendix. GF closed rings by Driss
  Bennis}
\stepcounter{section}

\noindent
Recall that a ring $R$ is said to be left GF-closed if in every exact
sequence of $R$-modules $0 \to G \to M \to G' \to 0$ with $G$ and $G'$
Gorenstein flat, also the middle module is Gorenstein flat. No example
is known of a ring that is not GF-closed.

From \lemcite[2.5]{DBn09} it is known that the middle module is
Gorenstein flat in every exact sequence $0 \to G \to M \to F \to 0$
with $F$ flat and $G$ Gorenstein flat.

\begin{thm}
  For a ring $R$ the following conditions are equivalent.
  \begin{eqc}
  \item $R$ is left GF-closed.
  \item In every exact sequence of $R$-modules $0 \to F \to M \to G
    \to 0$ with $F$ flat and $G$ Gorenstein flat, the module $M$ is
    Gorenstein flat.
  \end{eqc}
\end{thm}

\begin{prf*}
  The implication \proofofimp[]{i}{ii} is trivial. To prove that
  \eqclbl{ii} implies \eqclbl{i}, let
  \begin{equation*}
    0 \lra A \lra B \lra C \lra 0
  \end{equation*}
  be an exact sequence of $R$-modules where $A$ and $C$ are
  G-flat. The goal is to prove that $B$ is G-flat. For every injective
  $\Rop$-module $I$ one evidently has $\Tor{>0}{I}{B}=0$.  Thus,
  \pgref{g-flat}.(1) holds and we proceed to verify
  \pgref{g-flat}.(2).

  Since $A$ is G-flat, there is an exact sequence
  \begin{equation*}
    0 \lra A \lra F' \lra G' \lra 0\:
  \end{equation*}
  with $F'$ flat and $G'$ G-flat. From the push-out diagram
  \begin{equation*}
    \xymatrix@=1.5pc{
      & 0 \ar[d] & 0 \ar[d]\\
      0 \ar[r] & A \ar[r] \ar[d] & B \ar[r] \ar[d]
      & C \ar@{=}[d] \ar[r] & 0\\
      0 \ar[r] & F' \ar[d] \ar[r] & X \ar[r] \ar[d] & C \ar[r] & 0\\
      & G' \ar@{=}[r] \ar[d] & G' \ar[d] \\
      & 0& 0
    }
  \end{equation*}
  one gets by \eqclbl{ii} that $X$ is G-flat, and so there exists an
  exact sequence
  \begin{equation*}
    0 \lra X \lra F \lra G \lra 0
  \end{equation*}
  with $F$ flat and $G$ G-flat. Consider the second non-zero row in
  the push-out diagram
  \begin{equation*}
    \xymatrix@=1.5pc{
      && 0 \ar[d] & 0 \ar[d]\\
      0 \ar[r] & B \ar[r] \ar@{=}[d] & X \ar[r] \ar[d]
      & G' \ar[d] \ar[r] & 0\\
      0 \ar[r] & B \ar[r] & F \ar[r] \ar[d] & B' \ar[r] \ar[d] & 0\\
      && G \ar@{=}[r] \ar[d] & G \ar[d] \\
      && 0& 0\xycomma[.]
    }
  \end{equation*}
  The module $B'$ has the same properties as $B$; indeed $G$ and $G'$
  are G-flat, and hence one has $\Tor{>0}{I}{B'}=0$ for every
  injective $\Rop$-module $I$. Thus, by repeating the above process
  one can construct the complex in \pgref{g-flat}.(2).
\end{prf*}

\bibliographystyle{amsplain-nodash}


\def\soft#1{\leavevmode\setbox0=\hbox{h}\dimen7=\ht0\advance \dimen7
  by-1ex\relax\if t#1\relax\rlap{\raise.6\dimen7
  \hbox{\kern.3ex\char'47}}#1\relax\else\if T#1\relax
  \rlap{\raise.5\dimen7\hbox{\kern1.3ex\char'47}}#1\relax \else\if
  d#1\relax\rlap{\raise.5\dimen7\hbox{\kern.9ex \char'47}}#1\relax\else\if
  D#1\relax\rlap{\raise.5\dimen7 \hbox{\kern1.4ex\char'47}}#1\relax\else\if
  l#1\relax \rlap{\raise.5\dimen7\hbox{\kern.4ex\char'47}}#1\relax \else\if
  L#1\relax\rlap{\raise.5\dimen7\hbox{\kern.7ex
  \char'47}}#1\relax\else\message{accent \string\soft \space #1 not
  defined!}#1\relax\fi\fi\fi\fi\fi\fi}
  \providecommand{\MR}[1]{\mbox{\href{http://www.ams.org/mathscinet-getitem?mr=#1}{#1}}}
  \renewcommand{\MR}[1]{\mbox{\href{http://www.ams.org/mathscinet-getitem?mr=#1}{#1}}}
  \providecommand{\arxiv}[2][AC]{\mbox{\href{http://arxiv.org/abs/#2}{\sf
  arXiv:#2 [math.#1]}}} \def\cprime{$'$}
\providecommand{\bysame}{\leavevmode\hbox to3em{\hrulefill}\thinspace}
\providecommand{\MR}{\relax\ifhmode\unskip\space\fi MR }
\providecommand{\MRhref}[2]{%
  \href{http://www.ams.org/mathscinet-getitem?mr=#1}{#2}
}
\providecommand{\href}[2]{#2}

\end{document}